\documentclass[12pt]{article}
\usepackage{graphics, graphicx}     

\usepackage{cite}
\usepackage{amsmath}
\usepackage{amssymb}
\usepackage{amsthm}
\theoremstyle{plain}
\newtheorem{Theorem}{Theorem}[section] %
\newtheorem{Lemma}{Lemma}[section]
\newtheorem{Proposition}{Proposition}[section]

\theoremstyle{definition}
\newtheorem{Remark}{Remark}[section]
\newtheorem{Corollary}{Corollary}[section]

\theoremstyle{definition}

\newenvironment{Proof} 
{\par\noindent{\it Proof of}} 
{\hfill$\vspace{5mm}\scriptstyle\blacksquare$} 

\numberwithin{equation}{section} 
\numberwithin{figure}{section} 
\numberwithin{table}{section} 

\begin{document}

\setcounter{page}{1}

\markboth{M.I. Isaev, R.G. Novikov}{Reconstruction of a potential from the impedance boundary map}

\title{Reconstruction of a potential from the impedance boundary map}
\date{}
\author{ { M.I. Isaev and R.G. Novikov}}

\maketitle

\begin{abstract}
	We consider the inverse boundary value problem for the  Schr\"odinger equation
	at fixed energy with boundary measurements represented as the impedance boundary map
	(or Robin-to-Robin map). 
	We give formulas and equations for finding (generalized) scattering data for 
	the  aforementioned equation from boundary measurements in this impedance representation. 
	Combining these results with results of
	the inverse scattering theory we obtain efficient methods for reconstructing potential 
	from the impedance boundary map. To our knowledge, results of the present work are new 
	already for the case of Neumann-to-Dirichlet map.
\end{abstract}

\section{Introduction}

We consider the equation
\begin{equation}\label{eq} 
	-\Delta \psi  + v(x)\psi = E\psi, \ \  x \in  D, \, E\in\mathbb{R},
\end{equation}
where
\begin{equation}\label{eq_c}
	\begin{aligned}
		 D &\text{ is an open bounded domain in } \mathbb{R}^d,\ d\geq 2, \ \\ &\text{with } \partial D \in C^2,
 	\end{aligned}
\end{equation}
\begin{equation}\label{eq_c2}
 v \in \mathbb{L}^{\infty}(D), \ \ v = \bar{v}.
\end{equation}
Equation \eqref{eq} can be considered as the stationary Schr\"odinger equation of 
quantum mechanics at fixed energy $E$. 
Equation \eqref{eq} at fixed $E$ arises also in acoustics and 
electrodynamics.

Following  \cite{GM2009}, \cite{IN2012},
we consider the impedance boundary map $\hat{M}_\alpha = \hat{M}_{\alpha,v}(E)$ defined by
\begin{equation}\label{def_M}
	\hat{M}_\alpha [\psi]_\alpha = [\psi]_{\alpha-\pi/2}
\end{equation}
 for all sufficiently regular solutions $\psi$ of  equation (\ref{eq}) in $\bar{D} = D \cup \partial D$, where
\begin{equation}\label{[psi]_alpha}
	[\psi]_\alpha = [\psi(x)]_\alpha = \cos\alpha  \, \psi(x)  -  \sin\alpha\, \frac{\partial \psi}{\partial \nu}|_{\partial D}(x),\ \  
	x \in \partial D, \ \alpha \in \mathbb{R}
\end{equation}
and $\nu$ is the outward normal to $\partial D$. Under assumptions \eqref{eq_c}, \eqref{eq_c2},
in Lemma 3.2 of \cite{IN2012} it was shown 
that there is not more than a countable number of $\alpha\in\mathbb{R}$ such that 
$E$ is an eigenvalue for the operator $-\Delta + v$ in $D$
with the boundary condition 
\begin{equation}\label{bound_cond}
 \cos\alpha  \, \psi|_{\partial D}  -  \sin\alpha\, \frac{\partial \psi}{\partial \nu}|_{\partial D}  = 0.
\end{equation} 
Therefore, for any fixed $E$ we can assume that  for some fixed  $\alpha\in\mathbb{R}$
\begin{equation}\label{correct_M}
	\begin{aligned}
		\text{ $E$ is not an eigenvalue for the operator $-\Delta + v$ in $D$} \\
		\text{with boundary condition (\ref{bound_cond})}
	\end{aligned}
\end{equation}
and, as a corollary, $\hat{M}_\alpha$ can be defined correctly.

We consider $\hat{M}_\alpha = \hat{M}_{\alpha,v}(E)$
as an operator representation of all possible boundary measurements for the physical model
described by \eqref{eq}. We recall that the impedance boundary map $\hat{M}_\alpha$ is reduced to the 
Dirichlet-to-Neumann(DtN) map if $\alpha=0$ and is reduced to the 
Neumann-to-Dirichlet(NtD) map if $\alpha=\pi/2$. 
The map $\hat{M}_\alpha$ can be called also as the Robin-to-Robin map.

As in \cite{IN2012}, we consider the following inverse boundary value problem for equation (\ref{eq}).

{\bf Problem 1.1.}
 Given $\hat{M}_\alpha$ for some fixed $E$ and $\alpha$, find $v$. 

This problem can be considered as the Gel'fand inverse boundary value problem for the Schr\"odinger equation at fixed energy
(see \cite{Gelfand1954}, \cite{Novikov1988}). Note that in the initial Gel'fand formulation energy $E$
was not yet fixed and boundary measurements were considered as an operator relating $\psi|_{\partial D}$ and 
$\frac{\partial \psi}{\partial \nu}|_{\partial D}$ for $\psi$ satisfying \eqref{eq}.

Problem 1.1 for $E=0$ can be considered also as a generalization of the Calderon problem of the electrical impedance tomography (see \cite{Calderon1980}, \cite{Novikov1988}).

Note also that Problem 1.1 can be considered as an example of ill-posed problem: 
see \cite{BK2012}, \cite{LR1986}  for an introduction to this theory.

Problem 1.1 includes, in particular, the following questions: (a) uniqueness, (b) reconstruction,
(c) stability.

Global uniqueness theorems and global reconstruction methods for Problem 1.1 with $\alpha =0$ (i.e. for the DtN case)
were given for the first time in \cite{Novikov1988} in dimension $d \geq 3$  and in \cite{Buckhgeim2008} in dimension $d=2$.

Global stability estimates for Problem 1.1 with $\alpha =0$ were given for the first time in \cite{Alessandrini1988} 
in dimension $d \geq 3$ and in \cite{NS2010} in dimension $d=2$.  
A principal improvement of the result of \cite{Alessandrini1988} was given recently in \cite{Novikov2011} (for $E=0$).
Due to 
  \cite{Mandache2001} these logarithmic 
stability results are optimal (up to the value of the exponent). 
An extention of the instability estimates of  \cite{Mandache2001}  to the case of  non-zero energy as well as to the case 
 of Dirichlet-to-Neumann map given on the energy intervals was obtained in \cite{IsaevDtN}. 
 An extention of stability estimates of \cite{Novikov2011} to the energy dependent case was given recently in \cite{IN2012++}.
 Instability estimates complementing stability results of \cite{IN2012++} were obtained in \cite{I2012}.

Note also that for the Calderon problem (of the electrical impedance tomography) in its initial formulation the global uniqueness 
was firstly proved in \cite{SU1987} for $d\geq 3$ and in \cite{Nachman1996} for $d=2$. In addition, for the case of
piecewise constant or piecewise real analytic conductivity the first uniqueness results for the Calderon problem in dimension $d\geq 2$
were given in \cite{Druskin1982}, \cite{KV1985}. 
Lipschitz stability estimate for the case of
piecewise constant conductivity was proved in \cite{AV2005} and additional studies in this direction
were fulfilled in \cite{Rondi2006}.
 
 It should be noted that in  most of previous works on inverse boundary value problems for equation (\ref{eq}) at fixed $E$
 it was assumed in one way or another that 
 $E$ is not a Dirichlet eigenvalue for the operator $-\Delta + v$ in $D$, see \cite{Alessandrini1988}, \cite{Mandache2001}, \cite{Novikov1988}, \cite{Novikov2011}-\cite{S2011}. Nevertheless, the results of \cite{Buckhgeim2008} can be
 considered as global uniqueness  and reconstruction results for Problem 1.1 in dimension $d=2$ with general $\alpha$.
 
 Global stability estimates for Problem 1.1 in dimension $d\geq 2$ with general $\alpha$ 
 were recently given in \cite{IN2012}.

In the present work we give formulas and equations for finding (generalized) scattering data 
from the impedance boundary map $\hat{M}_\alpha$ with general $\alpha$. Combining these results with results of 
\cite{Grinevich1988}, \cite{Henkin1987}, \cite{Nachman1996}, \cite{Novikov1992}-\cite{Novikov1999}, \cite{Novikov2005bar}-\cite{Novikov2009},
we obtain efficient reconstruction methods for Problem 1.1 in multidimensions with general $\alpha$. To
our knowledge these results are new already for the NtD case.

In particular, in the present work we give first mathematically justified approach for reconstructing coefficient
$v$ from boundary measurements for \eqref{eq} via inverse scattering without the assumption that $E$ is not a Dirichlet 
eigenvalue for $-\Delta + v$ in $D$. In addition, numerical efficiency of related inverse scattering 
techniques was shown in \cite{ABR2008}, \cite{BBMRS2000}, \cite{BAR2009}, \cite{BSRZ2012}; see also \cite{BKM2011}.

Definitions of (generalized) scattering data are recalled in Section 2. Our main results are presented in Section 3. 
Proofs of these results are given in Sections 4, 5 and 6.


\section{Scattering data}
Consider the Schr\"odinger equation
\begin{equation}\label{2.1}
-\Delta\psi+v(x)\psi=E\psi,\ \ x\in \mathbb{R}^d, \ \ d\ge 2
\end{equation}
where
\begin{equation}\label{2.2}
 (1+|x|)^{{d+\varepsilon}}v(x)\in \mathbb{L}^{\infty}(\mathbb{R}^d) \text{ (as a function of $x$)},\ \ {\rm for\ some}\ \ \varepsilon>0.
\end{equation}
 For equation \eqref{2.1} we consider the functions $\psi^+$ and $f$ of the
classical scattering theory and the Faddeev functions $\psi$, $h$,
$\psi_{\gamma}$, $h_{\gamma}$ (see, for example, \cite{Berezin1991}, \cite{E2011}, 
 \cite{Faddeev1965}, \cite{Faddeev1974}, \cite{Grinevich2000},   \cite{Henkin1987}, \cite{Newton1989}, 
  \cite{Novikov1992}).

The functions $\psi^+$ and $f$ can be defined as follows:

\begin{equation}\label{2.3}
	\psi^+(x,k) = e^{ikx} + \int\limits_{\mathbb{R}^d} G^+(x-y,k) v(y)\psi^+(y,k)dy,
\end{equation}
\begin{equation}\label{2.4}
\begin{aligned}
	G^+(x,k) = -\left(\frac{1}{2\pi}\right)^d \int\limits_{\mathbb{R}^d}\frac{e^{i\xi x}}{\xi^2 - k^2 -i0} \, d\xi,\\
	x,k\in\mathbb{R}^d,\ \  k^2>0,
\end{aligned}
\end{equation}

where (\ref{2.3}) at fixed $k$  is  considered as an
equation for $\psi^+$ in $\mathbb{L}^{\infty}(\mathbb{R}^d)$;
\begin{equation}\label{2.5}
	\begin{aligned}
	f(k,l) = \left(\frac{1}{2\pi}\right)^d \int\limits_{\mathbb{R}^d} e^{-ilx} \psi^+(x,k)v(x)dx, 
	\\ k,l\in \mathbb{R}^d, \ k^2>0.
	\end{aligned}
\end{equation} 
 In addition: $\psi^+(x,k)$ satisfies (\ref{2.3}) for
$E=k^2$ and describes scattering of the plane waves $e^{ikx}$; $f(k,l)$,
$k^2=l^2$, is the scattering amplitude for equation (\ref{2.1}) for $E=k^2$.
Equation (\ref{2.3}) is the Lippman-Schwinger integral equation.

The functions $\psi$ and $h$ can be defined as follows:
	\begin{equation}\label{eq_psi}
			\psi(x,k) = e^{ikx} + \int\limits_{\mathbb{R}^d} G(x-y,k) v(y)\psi(y,k)dy,
	\end{equation}
	\begin{equation}\label{G(x,k)}
	\begin{aligned}
			G(x,k) = -\left(\frac{1}{2\pi}\right)^d \int\limits_{\mathbb{R}^d}\frac{e^{i\xi x} d\xi}{\xi^2 + 2k\xi}\, e^{ikx}, \\
			x\in\mathbb{R}^d, \ \ k\in\mathbb{C}^d, \ \ \mbox{Im}\, k \neq 0, 
	\end{aligned}
	\end{equation}

where (\ref{eq_psi}) at fixed $k$  is considered as
an equation for $\psi=e^{ikx}\mu(x,k)$, $\mu \in \mathbb{L}^{\infty}(\mathbb{R}^d)$; 

\begin{equation}\label{h(k,l)}
		\begin{aligned}
				h(k,l) = \left(\frac{1}{2\pi}\right)^d \int\limits_{\mathbb{R}^d} e^{-ilx} \psi(x,k)v(x)dx, \\ 
				k,l\in \mathbb{C}^d, 
				\  \mbox{Im}\, k = \mbox{Im}\, l \neq 0.
		\end{aligned}
	\end{equation} 
 In addition, $\psi(x,k)$ satisfies
(\ref{2.1}) for $E=k^2$, and $\psi$, $G$ and $h$ are (nonanalytic)
continuations of $\psi^+$, $G^+$ and $f$ to the complex domain. In particular,
 $h(k,l)$ for $k^2=l^2$ can be considered as the "scattering" amplitude
in the complex domain for equation \eqref{2.1} for $E=k^2$.
The functions $\psi_{\gamma}$ and $h_{\gamma}$ are defined as follows:
\begin{equation}\label{2.9}
\begin{aligned}
\psi_{\gamma}(x,k)=\psi(x,k+i0\gamma),\ \ \text{ }\  \ h_{\gamma}(k,l)=h(k+i0\gamma,
l+i0\gamma),\\
x,k,l,\gamma\in\mathbb{R}^d, \  \gamma^2=1.
\end{aligned}
\end{equation}
 We recall also that
\begin{equation}
\begin{aligned}
\psi^+(x,k)=\psi_{k/|k|}(x,k),\ \ \text{ }\ f(k,l)=h_{k/|k|}(k,l),\\
x,k,l\in\mathbb{R}^d, \ |k|>0.
\end{aligned}
\end{equation}

We consider $f(k,l)$ and $h_{\gamma}(k,l)$, where $k,l,\gamma\in\mathbb{R}^d$,
$k^2=l^2=E$, $\gamma^2=1$, and $h(k,l)$, where $k,l\in\mathbb{C}^d$,
$\mbox{Im}\,k=\mbox{Im}\,l\neq 0$, $k^2=l^2=E$, as scattering data $S_E$ for equation (\ref{2.1}) at
fixed $E\in ( 0, +\infty )$. We consider $h(k,l)$, where $k,l\in\mathbb{C}^d$,
$\mbox{Im}\,k=\mbox{Im}\,l\neq 0$, $k^2=l^2=E$, as scattering data $S_E$ for equation (\ref{2.1}) at
fixed $E\in (-\infty,0]$.

We consider also the sets  ${\cal E}$, ${\cal E}_{\gamma}$,
${\cal E}^+$ defined as follows:
\begin{subequations}\label{2.11}
\begin{equation}\label{2.11a}
{\cal E}=\left\{
\begin{array}{c}
\zeta\in\mathbb{C}^d\setminus \mathbb{R}^d\ :\ \ {\rm equation \ (\ref{eq_psi})} \ \ {\rm for}\ \
k=\zeta\ \ {\rm is\ not}
\\ {\rm uniquely\ solvable\ for}\ \ \psi=e^{ikx}\mu\ \ {\rm with}\ \
\mu\in \mathbb{L}^{\infty}(\mathbb{R}^d)
\end{array}
\right\},
\end{equation}

\begin{equation}\label{2.11b}
\begin{aligned}
{\cal E}_\gamma=\left\{
\begin{array}{c}
\zeta\in\mathbb{R}^d\setminus \{0\}\ :\ \ {\rm equation \ (\ref{eq_psi})} \ \ {\rm for}\ \
k=\zeta+i0\gamma\ \\ {\rm is\ not} \
 {\rm uniquely\ solvable\ for}\ \ \psi=\mathbb{L}^{\infty}(\mathbb{R}^d)\
\end{array}
\right\}, \\ \gamma\in \mathbb{S}^{d-1},
\end{aligned}
\end{equation}

\begin{equation}\label{2.11c}
\begin{aligned}
{\cal E}^+=\left\{
\begin{array}{c}
\zeta\in\mathbb{R}^d\setminus \{0\}\ :\ \ {\rm equation \ (\ref{eq_psi})} \ \ {\rm for}\ \
k=\zeta\ \ \ {\rm is\ not}
\\ {\rm uniquely\ solvable\ for}\ \ \psi=\mathbb{L}^{\infty}(\mathbb{R}^d)\
\end{array}
\right\}.
\end{aligned}
\end{equation}
\end{subequations}

In addition, ${\cal E}^+$ is a well-known set of the classical scattering theory
for equation (\ref{2.1}) and ${\cal E}^+=\emptyset$ for real-valued $v$ satisfying
(\ref{2.2})
(see, for example, \cite{Berezin1991}, \cite{Newton1989}). 
Note also that $\mathcal{E}^+$ is spherically symmetric.
The sets ${\cal E}$, ${\cal E}_{\gamma}$ were
considered
for the first time in \cite{Faddeev1965}, \cite{Faddeev1974}. Concerning the properties of
${\cal E}$
and ${\cal E}_{\gamma}$, see \cite{Faddeev1974}, \cite{GN2012}, \cite{Henkin1987}, \cite{LN1987}, \cite{Newton1989}, 
 \cite{Nachman1996}, \cite{Novikov1996}, \cite{Weder1991}.

We consider also the functions $R$, $R_{\gamma}$,
$R^+$ defined as follows:
\begin{equation}\label{R(x,y,k)}
\begin{aligned}
R(x,y,k)=G(x-y,k)+\int\limits_{\mathbb{R}^d}G(x-z,k)v(z)R(z,y,k)dz,\\
x,y\in\mathbb{R}^d,\  k\in\mathbb{C}^d, \ \mbox{Im}\,k \neq 0,
\end{aligned}
\end{equation}
where $G$ is defined by (\ref{G(x,k)}) and formula \eqref{R(x,y,k)} at
fixed $y$, $k$ is considered as an equation for
\begin{equation}\label{2.13}
R(x,y,k)=e^{ik(x-y)}r(x,y,k),
\end{equation}
where $r$ is sought with the properties
\begin{subequations}\label{2.14}
\begin{equation}
r(\cdot,y,k)\ \ {\rm is\ continuous\ on}\ \ \mathbb{R}^d\setminus \{y\}
\end{equation}
\begin{equation}
r(x,y,k)\to 0\ \ {\rm as}\ \ |x|\to\infty,
\end{equation}
\begin{equation}
\begin{aligned}
&r(x,y,k)=O(|x-y|^{2-d})\ \ {\rm as}\ \ x\to y\ \ {\rm for}\ \ d\ge 3,\\
&r(x,y,k)=O(|\ln\,|x-y||)\ \ {\rm as}\ \ x\to y\ \ {\rm for}\ \ d=2;
\end{aligned}
\end{equation}
\end{subequations}
\begin{equation}\label{2.15}
\begin{aligned}
	R_{\gamma}(x,y,k)&=R(x,y,k+i0\gamma),\\
		&x,y\in\mathbb{R}^d, \ k\in\mathbb{R}^d\setminus \{0\},\  \gamma\in\mathbb{S}^{d-1};
\end{aligned}
\end{equation}
\begin{equation}\label{2.16}
\begin{aligned}
R^+(x,y,k)=R_{k/|k|}(&x,y,k),\\
 &x,y\in\mathbb{R}^d,\  k\in\mathbb{R}^d\setminus \{0\}.
\end{aligned}
\end{equation}
 In addition, the functions $R(x,y,k)$, $R_{\gamma}(x,y,k)$
and $R^+(x,y,k)$ (for their domains of definition in $k$ and $\gamma$)
satisfy the following equations:
\begin{equation}\label{2.17}
\begin{aligned}
	(\Delta_x+E-v(x))R(x,y,k)=\delta(x-y),\\
	(\Delta_y+E-v(y))R(x,y,k)=\delta(x-y),\\  x,y\in\mathbb{R}^d, \ E=k^2.
\end{aligned}
\end{equation}
The function $R^+(x,y,k)$ (defined by means of (\ref{R(x,y,k)}) for $k\in\mathbb{R}^d\setminus \{0\}$
with $G$ replaced by $G^+$ of (\ref{2.4})) is well-known in the scattering theory
for equations (\ref{2.1}), \eqref{2.17} (see, for example, \cite{Berezanskii1958}). In particular, this function  describes
scattering of the spherical waves $G^+(x-y,k)$ generated by a source at $y$.
In addition $R^+(x,y,k)$ is a radial function in $k$, i.e.
\begin{equation}\label{2.18}
	R^+(x,y,k) = \mbox{R}^+(x,y,|k|), \ \ \ \ x,y\in \mathbb{R}^d, \ k\in \mathbb{R}^d \setminus \{0\}.
\end{equation}
Apparently, the functions $R$ and $R_{\gamma}$ were considered for the
first time in \cite{Novikov1996}.

In addition, under the assumption (\ref{2.2}): equation (\ref{R(x,y,k)}) at fixed $y$ and $k$
is uniquely solvable for $R$ with the properties (\ref{2.13}), (\ref{2.14}) if and only if
$k\in\mathbb{C}^d\setminus (\setminus \mathbb{R}^d\cup {\cal E})$; equation (\ref{R(x,y,k)}) with $k=\zeta+i0\gamma$,
$\zeta\in\mathbb{R}^d\setminus \{0\}$, $\gamma\in\mathbb{S}^{d-1}$, at fixed $y$, $\zeta$ and $\gamma$
is uniquely solvable for $R_{\gamma}$ if and only if $\zeta\in\mathbb{R}^d\setminus
(\{0\}\cup {\cal E}_{\gamma})$; equation (\ref{R(x,y,k)}) with $k=\zeta+i0{\zeta/|\zeta|}$,
$\zeta\in\mathbb{R}^d\setminus {0}$, at fixed $y$ and $\zeta$ is uniquely solvable for
$R^+$ if and only if $\zeta\in\mathbb{R}^d\setminus (\{0\}\cup {\cal E}^+)$.


\section{Main results}

 Let  $v$ and $v^0$   satisfy  (\ref{eq_c2}), (\ref{correct_M}) for some fixed $E$ and $\alpha$.
 Let $M_{\alpha,v}(x,y,E)$, $M_{\alpha,v^0}(x,y,E)$, $x,y\in \partial D$, 
 denote the Schwartz kernels of the impedance boundary maps
 $\hat{M}_{\alpha,v}$, $\hat{M}_{\alpha,v^0}$, for potentials $v$ and $v^0$, respectively, where 
 $\hat{M}_{\alpha,v}$, $\hat{M}_{\alpha,v^0}$
are considered as linear integral operators. In addition, we consider 
$v^0$ as some known background potential.

Let $h$, $\psi$, $f$, $\psi^{+}$, $h_\gamma$, $\psi_\gamma$, $\mathcal{E}$, $\mathcal{E}^+$, $\mathcal{E}_\gamma$
and $h^0$, $\psi^0$, $f^0$, $\psi^{+,0}$, $h_\gamma^0$, $\psi_\gamma^0$, $\mathcal{E}^0$, $\mathcal{E}^{+,0}$, $\mathcal{E}_\gamma^0$
denote the functions and sets of \eqref{2.3}, \eqref{2.5}, \eqref{eq_psi}, \eqref{h(k,l)}, \eqref{2.9}, \eqref{2.11}
for  potentials $v$ and $v^0$, respectively. Here and bellow in this section we always assume that 
$v\equiv 0$, $v^0\equiv 0$ on $\mathbb{R}^d \setminus D$.

\begin{Theorem}\label{main}
Let $D$ satisfy \eqref{eq_c} and potentials $v$, $v^0$ satisfy (\ref{eq_c2}), (\ref{correct_M})
for some fixed $E$ and  $\alpha$. Then:
	\begin{equation}\label{eq_h}
		\begin{aligned}
		h(k,l) -  h^{0}(k,l) &=\\= \left(\frac{1}{2\pi}\right)^d 
		\int\limits_{\partial D}\int\limits_{\partial D} 
		&[\psi^{0}(x,-l)]_\alpha \left(M_{\alpha,v} - M_{\alpha,v_0}\right)(x,y,E) [\psi(y,k)]_\alpha dx \, dy,\\
			&\ \  \ \text{ $k,l\in \mathbb{C}^d\setminus (\mathcal{E} \cup \mathcal{E}^0)$, 
			$k^2 = l^2 = E$, $\mbox{Im}\, k = \mbox{Im}\, l\neq 0$, }
		\end{aligned}
	\end{equation}
		
	\begin{equation}\label{eq_psi_alpha}
		\begin{aligned}
		\	[\psi(x,k)]_\alpha = 
		[\psi^0(x,k)]_\alpha + 
		\int\limits_{\partial D} A_\alpha(x,y,k) [\psi(y,k)]_\alpha dy,\\
			\text{ $x \in \partial D$, $k\in \mathbb{C}^d\setminus (\mathcal{E} \cup \mathcal{E}^0)$, $\mbox{Im}\,k\neq 0$, $k^2=E$}
		\end{aligned}
	\end{equation}
	where
	\begin{equation}\label{A_alpha}
		A_\alpha(x,y,k) = \lim_{\varepsilon \rightarrow +0}
		\int\limits_{\partial D} D_{\alpha,\varepsilon} R^0(x,\xi,k) \left(M_{\alpha,v} - M_{\alpha,v^0}\right)(\xi,y,E) d\xi,	
	\end{equation}
	\begin{equation}\label{DalphaG}
		\begin{aligned}
		D_{\alpha,\varepsilon}R^0(x,\xi,k) =  [[R^0(x+\varepsilon \nu_x,\xi,k)]_{\xi,\alpha}]_{x,\alpha} &= \\ =\Bigg(
			\cos^2\alpha  
			- \sin \alpha \cos \alpha \left( \frac{\partial}{\partial \nu_x} +  \frac{\partial}{\partial \nu_\xi}\right) 
			+ \sin^2\alpha \frac{\partial^2}{\partial \nu_x \partial \nu_\xi}
		\Bigg)  &R^0(x+\varepsilon \nu_x, \xi,k),\\
			&x,\xi,y \in \partial D,
		\end{aligned}
	\end{equation}
	where $R^0$ denotes the Green function of \eqref{R(x,y,k)} for potential $v^0$,
	$\nu_x$ is the outward normal to $\partial D$ at $x$.
	In addition,
	formulas completely similar to \eqref{eq_h} - \eqref{DalphaG} are also valid for 
	the classical scattering functions
	$f$, $\psi^+$, $f^0$, $\psi^{+,0}$ and sets $\mathcal{E}^+$, $\mathcal{E}^{+,0}$  of 
	\eqref{2.3}, \eqref{2.5}, \eqref{2.11c}
	for $v$ and $v^0$, respectively, but with $R^{+,0}$ in place of $R^0$ in \eqref{A_alpha}, \eqref{DalphaG}, where
	$R^{+,0}$ denotes the Green function of \eqref{2.16} for potential $v^0$.
\end{Theorem}

Theorem 3.1 is proved in Section 4.

Note that formula of the type \eqref{eq_h} for $h_\gamma$ is not completely similar to \eqref{eq_h}: 
see formula \eqref{eq_h_gamma} given below. In 
this formula \eqref{eq_h_gamma}, in addition to expected $\psi_\gamma(x,k)$,
 we use also $\psi_\gamma(x,k,l)$ defined as follows:

\begin{equation}\label{psi(x,k,l)}
	\begin{aligned}
		\psi_\gamma(x,k,l) = e^{ilx} + \int\limits_{\mathbb{R}^d} G_\gamma(x-y,k) v(y)\psi_\gamma(y,k,l)dy,\\
	 G_\gamma(x,k) = G(x,k+i0\gamma),
		\\
		\gamma \in \mathbb{S}^{d-1}, \ x,k,l\in \mathbb{R}^d, k^2=l^2>0,
	\end{aligned}
\end{equation}
where \eqref{psi(x,k,l)} at fixed $\gamma$, $k$, $l$ is considered as an equation for  
$\psi_\gamma(\cdot,k,l)$
in $\mathbb{L}^\infty(\mathbb{R}^d)$, 
$G$ is defined by \eqref{G(x,k)}.
\begin{Proposition}
Let the asssumptions of Theorem \ref{main} hold. Let $\psi_\gamma(x,k)$
correspond to $v$ according to \eqref{2.9} and $\psi^0_{-\gamma}(\cdot,k,l)$
correspond to $v^0$
according to \eqref{psi(x,k,l)}. Then
\begin{equation}\label{eq_h_gamma}
		\begin{aligned}
		h_\gamma(k,l) -  h^{0}_\gamma(k,l) &=\\= \left(\frac{1}{2\pi}\right)^d 
		\int\limits_{\partial D}\int\limits_{\partial D} 
		&[\psi_{-\gamma}^{0}(x,-k,-l)]_\alpha \left(M_{\alpha,v} - M_{\alpha,v^0}\right)(x,y,E) [\psi_\gamma(y,k)]_\alpha dx \, dy,\\
			&\ \  \ \text{$\gamma \in \mathbb{S}^{d-1}$, $k\in \mathbb{R}^d \setminus (\{0\}\cup 
			\mathcal{E}_\gamma \cup \mathcal{E}_\gamma^0)$, $l\in\mathbb{R}^d$, 
			$k^2 = l^2=E$. }
		\end{aligned}
	\end{equation}
	In addition, formulas completely similar to \eqref{eq_psi_alpha} - \eqref{DalphaG} are also valid for the functions
	$\psi_\gamma(x,k)$, $\psi^0_\gamma(x,k)$ and sets $\mathcal{E}_\gamma$, $\mathcal{E}_\gamma^0$
	of \eqref{2.9}, \eqref{2.11b} for $v$ and $v^0$, respectively, but with $R_\gamma^0$ in place of $R^0$ in
	\eqref{A_alpha}, \eqref{DalphaG}, where
	$R_\gamma^{0}$ denotes the Green function of \eqref{2.15} for potential $v^0$.
\end{Proposition}

Proposition 3.1 is proved in Section 4.

Note that \eqref{eq_psi_alpha} is considered as a linear integral equation for finding $[\psi(x,k)]_\alpha$, $x\in \partial D$,
at fixed $k$, from $\hat{M}_{\alpha,v}-\hat{M}_{\alpha,v^0}$ and
$[\psi^0(x,k)]_\alpha$, whereas (\ref{eq_h})
is considered as an explicit formula for finding $h$ from $h^0$, $\hat{M}_{\alpha,v} - \hat{M}_{\alpha,v^0}$, $[\psi^0(x,k)]_\alpha$
 and $[\psi(x,k)]_\alpha$. 
In addition, we use similar interpretation for similar formulas for $\psi^+$, $f$ and for $\psi_\gamma$, $h_\gamma$,
mentioned in Theorem 3.1 and Proposition 3.1.

Under the assumptions of Theorem \ref{main}, the following propositions are valid:

\begin{Proposition} Equation \eqref{eq_psi_alpha} for $[\psi(x,k)]_\alpha$ at fixed 
$k\in \mathbb{C}^d\setminus (\mathbb{R}^d\cup \mathcal{E}^0)$ is 
a Fredholm linear integral equation of the second kind in the space of 
bounded functions on $\partial D$. In addition, the same is also valid for the equation
for $[\psi^+(x,k)]_\alpha$ at fixed $k\in \mathbb{R}^d\setminus (\{0\}\cup\mathcal{E}^{+,0})$,
mentioned in Theorem \ref{main}, and for the equation for $[\psi_\gamma(x,k)]_\alpha$ 
at fixed $\gamma \in \mathbb{S}^{d-1}$, $k\in \mathbb{R}^d\setminus (\{0\}\cup\mathcal{E}_\gamma^{0})$,
mentioned in Proposition 3.1. 
\end{Proposition}

Proposition 3.2 is proved in Section 4.

\begin{Proposition} 
For $k\in \mathbb{C}^d\setminus (\mathbb{R}^d\cup \mathcal{E}^0)$
equation
\eqref{eq_psi_alpha} is uniquely solvable in the space of bounded functions 
on  $\partial D$ if and only if $k \notin \mathcal{E}$. In addition, 
the aforementioned equations for  $[\psi^+(x,k)]_\alpha$, $k\in \mathbb{R}^d\setminus (\{0\}\cup\mathcal{E}^{+,0})$,
and $[\psi_\gamma(x,k)]_\alpha$, $\gamma\in \mathbb{S}^{d-1}$, $k\in \mathbb{R}^d\setminus (\{0\}\cup\mathcal{E}_\gamma^{0})$,  
are uniquely solvable in the space of bounded functions on $\partial D$ if and only if
$k\notin \mathcal{E}^+$ and $k\notin \mathcal{E}_\gamma$, respectively.
\end{Proposition}

Proposition 3.3 is proved in Section 5.

\begin{Proposition}
Let $\phi_\alpha(x,y)$ 
be the solution of the Dirichlet boundary value problem at fixed $y\in \partial D$, $\lambda \in \mathbb{C}$:
\begin{equation}\label{eq_phi}
\begin{aligned}
	&-\Delta_x \phi_\alpha(x,y) = \lambda \phi_\alpha(x,y), \ \ \text{ } \ \ x\in D, \\  
	&\phi_\alpha(x,y) = \left(M_{\alpha,v} - M_{\alpha,v^0}\right)(x,y,E), \ \ \ x\in \partial D,
\end{aligned}
\end{equation}
where we assume that $\lambda$ is not a Dirichlet eigenvalue for $-\Delta$ in $D$.
Then
\begin{equation}\label{3.8}
\begin{aligned}
		A_\alpha(x,y,k) &= \lim\limits_{\varepsilon \rightarrow +0} 
		\int\limits_{\partial D}  [R^0(x+\varepsilon\nu_x,\xi,k)]_{x,\alpha} [\phi_{\alpha}(\xi,y)]_{\xi,\alpha} d\xi -\\&-
 \sin \alpha  \int\limits_{D}  [R^0(x,\xi,k)]_{x,\alpha}(v^0(\xi)- E + \lambda) \phi_{\alpha}(\xi,y) d\xi, \ \ \ x,y\in \partial D,
\end{aligned}
\end{equation}
where 
\begin{equation}
		\ [R^0(x+\varepsilon\nu_x,\xi,k)]_{x,\alpha} = \left(\cos \alpha  
		- \sin\alpha \frac{\partial}{\partial \nu_{x}}\right) R^0(x+\varepsilon\nu_x,\xi,k), \ \ \ x\in \partial D,\ \xi \in \bar{D},
\end{equation}
\begin{equation}\label{3.10}
		\begin{aligned} 
			\ [\phi_\alpha(\xi,y)]_{\xi,\alpha} &= \left(\cos \alpha  
				- \sin\alpha \frac{\partial}{\partial \nu_{\xi}}\right)\phi_\alpha(\xi,y) =\\
				&=\cos \alpha\,\phi_\alpha(\xi,y)  - 
				\sin \alpha \left(\hat{\Phi}(\lambda)\phi_\alpha(\cdot,y)\right)(\xi), \ \ \ \xi,y \in \partial D,
	\end{aligned}
\end{equation}
where $A_\alpha$ is defined in \eqref{A_alpha}, 
$\hat{\Phi}(\lambda) = \hat{M}_{0,0}(\lambda)$ is the Dirichlet-to-Neumann map for \eqref{eq_phi}. In addition, formulas completely similar to \eqref{3.8} are also
valid for the kernels $A_\alpha^+$
(but with $R^+_0$ in place of $R_0$) and $A_{\alpha,\gamma}$
(but with $R_\gamma^{0}$ in place of $R^0$), arising in the equations for $[\psi^+]_\alpha$ and $[\psi_\gamma]_\alpha$,
mentioned in Theorem 3.1 and Proposition 3.1.
\end{Proposition}

Proposition 3.4 is proved in Section 4.

Note that, for the case when $\sin \alpha =0$, formula \eqref{3.8} coincides with \eqref{A_alpha}. However, for $\sin \alpha \neq 0$,
formula \eqref{3.8} does not contain $\partial^2 R^0/\partial \nu_x \partial \nu_\xi$ in contrast with \eqref{A_alpha}
and is more convenient than \eqref{3.8} in this sense.

Theorem 3.1, Propositions 3.1 - 3.4 and the reconstruction results from generalized
scattering data
 (see \cite{Grinevich2000}, \cite{Grinevich1988}, \cite{Henkin1987}, \cite{Novikov1992}-\cite{Novikov1999}, 
 \cite{Novikov2005bar}-\cite{Novikov2009}, 
 \cite{NS2012}) imply the following corollary: 
 
\begin{Corollary} To reconstruct a potential $v$ in the domain $D$ from its impedance boundary map $\hat{M}_{\alpha,v}(E)$  
at fixed $E$ and $\alpha$ one can use the following schema:

\begin{enumerate} 
\item $v^0\rightarrow  \{S_E^0\}, \{R^0\}, \{[\psi^0]_\alpha\}, \hat{M}_{\alpha,v^0}$ via direct problem methods,

\item $\{R^0\}, \hat{M}_{\alpha,v^0}, \hat{M}_{\alpha,v} \rightarrow \{A_\alpha\}$ as described in 
Theorem 3.1 and Propositions 3.1, 3.4,

\item $\{A_\alpha\}, \{[\psi^0]_\alpha\} \rightarrow \{[\psi]_\alpha\}$ as described
in 
Theorem 3.1 and Proposition 3.1,

\item $\{S_E^0\}, \{[\psi^0]_\alpha\}, \{[\psi]_\alpha\},
\hat{M}_{\alpha,v^0}, \hat{M}_{\alpha,v}
\rightarrow \{S_E\}$  as described 
in 
Theorem 3.1 and Proposition 3.1, 

\item $\{S_E\} \rightarrow v$ as described in \cite{Grinevich2000}, \cite{Grinevich1988}, \cite{Henkin1987}, \cite{Novikov1992}-\cite{Novikov1999}, \cite{Novikov2005bar}-\cite{Novikov2009}, 
 \cite{NS2012},
\end{enumerate}
where $\{S_E^0\}$ and
$\{S_E\}$
denote some appropriate
part of  
$h^0$,
$f^0$,
$h^0_\gamma$ and
$h$, $f$, $h_\gamma$, respectively,
$\{[\psi^0]_\alpha\}$ and
$\{[\psi]_\alpha\}$
denote some appropriate
part of  
$[\psi^0]_\alpha$,
$[\psi^{+,0}]_\alpha$,
$[\psi^0_\gamma]_\alpha$ and 
$[\psi]_\alpha$,
$[\psi^+]_\alpha$,
$[\psi_\gamma]_\alpha$, respectively,
$\{R^0\}$, 
$\{A_\alpha\}$ 
denote some appropriate
part of 
$R^0$, $R^{+,0}$,
$R^0_\gamma$,
$A_\alpha$, $A^+_\alpha$, $A_{\alpha,\gamma}$.
\end{Corollary}

\begin{Remark}
	For the case when $v^0 \equiv 0$, $\sin \alpha =0$, Theorem 3.1, Propositions 3.1 - 3.3
	and Corollary 3.1 (with available references at that time at step 5) were obtained in \cite{Novikov1988} 
	(see also \cite{Nachman1988}, \cite{Nachman1996}). Note that
	basic results of \cite{Novikov1988} were presented already in the survey given in \cite{Henkin1987}.
	For the case when $\sin \alpha =0$ Theorem 3.1, Propositions 3.1 - 3.3
	and Corollary 3.1 (with available references at that time at step 5) were obtained in \cite{Novikov2005}.
\end{Remark}

\begin{Remark}
	The results of Theorem 3.1, Propositions 3.1 - 3.4 and Corollary 3.1 remain valid
	for complex-valued $v$, $v^0$ and complex $E$, $\alpha$, under the condition that \eqref{correct_M} holds for
	both $v$ and $v^0$.
\end{Remark}

\begin{Remark}
	Under the assumptions of Theorem 3.1, the following formula holds:
	\begin{equation}\label{3.11}
		\hat{M}_{\alpha,v}(E) - \hat{M}_{\alpha,v^0}(E) = (D_\alpha \mbox{R}^{+,0}(E))^{-1} - (D_\alpha \mbox{R}^{+}(E))^{-1},
	\end{equation}
	\begin{equation}
		\begin{aligned}
		D_\alpha \mbox{R}^{+}(E) u(x) = \lim\limits_{\varepsilon \rightarrow +0}
		\int\limits_{\partial D} D_{\alpha,\varepsilon} \mbox{R}^{+}(x,y,\sqrt{E}) u(y) dy,\\
		D_\alpha \mbox{R}^{+,0}(E) u(x) = \lim\limits_{\varepsilon \rightarrow +0}
		\int\limits_{\partial D} D_{\alpha,\varepsilon} \mbox{R}^{+,0}(x,y,\sqrt{E}) u(y) dy, \\ x\in\partial D,
		\end{aligned}
	\end{equation}
	where $D_{\alpha, \varepsilon}$ is defined as in \eqref{DalphaG}, $\mbox{R}^{+}(x,y,\sqrt{E})$, $\mbox{R}^{+,0}(x,y,\sqrt{E})$, $\sqrt{E}>0$, 
	are the Green functions of \eqref{2.16} written as in \eqref{2.18} for potentials $v$, $v^0$, respectively, 
	$u$ is the test function. For the case when $\sin \alpha=0$, $v^0 \equiv 0$, $d\geq 3$, formula \eqref{3.11} was given in \cite{Nachman1988}.
	Using techniques developed in \cite{IN2012} and in the present work, we obtain \eqref{3.11} in the general case.
\end{Remark}

\section{Proofs of Theorem 3.1 and Propositions 3.1, 3.2,  3.4}

In this section we will use formulas and equations for impedance boundary map from \cite{IN2012}.
These results are presented in detail in Subsection 4.1. Proofs of Theorem 3.1 and
Propositions 3.1, 3.2, 3.4 are given in Subsections 4.2, 4.3.

\subsection{Preliminaries}
Let	$G_{\alpha,v}(x,y,E)$ be the Green function for the operator $\Delta - v + E$ in $D$
with the impedance
boundary condition \eqref{bound_cond} under assumptions \eqref{eq_c}, \eqref{eq_c2} and \eqref{correct_M}. 
We recall that (see formulas (3.12), (3.13) of \cite{IN2012}):
\begin{equation}\label{4.1}
	G_{\alpha,v}(x,y,E) = G_{\alpha,v}(y,x,E),  \ \ \ x,y \in \bar{D},
\end{equation}	
and, for $\sin \alpha \neq 0$,
\begin{equation}\label{4.2}
	M_{\alpha,v}(x,y,E) = \frac{1}{\sin^2\alpha} G_{\alpha,v}(x,y,E) - \frac{\cos \alpha}{\sin \alpha} \delta_{\partial D}(x-y), \ \ \ 
	x,y \in \partial D,
\end{equation}
where  $M_\alpha(x,y,E)$ and  $\delta_{\partial D}(x-y)$ denote 
the Schwartz kernels of the impedance boundary map $\hat{M}_{\alpha,v}(E)$
 and the identity operator $\hat{I}$ on $\partial D$, respectively, 
 where $\hat{M}_\alpha$ and $\hat{I}$ are considered as linear integral operators.
 
 We recall also that (see, for example, formula (3.16) of \cite{IN2012}):
 \begin{equation}\label{psiG+++}
 	\psi(x) = \frac{1}{\sin \alpha} \int\limits_{\partial D} (\cos\alpha\, \psi(\xi) - \sin\alpha \frac{\partial}{\partial \nu}\psi(\xi)) G_{\alpha,v}(x,\xi,E) d\xi,
 	\ \ \ x\in D,
 \end{equation} 
 for all sufficiently regular solutions $\psi$ of equation \eqref{eq} in $\bar{D}$ and $\sin \alpha \neq 0$.
 
 We will use the following properties of the Green function $G_\alpha(x,y,E)$:
 \begin{equation}\label{G_alpha_1}
	G_{\alpha,v}(x,y,E)\ \ {\rm is\ continuous\ in}\ \ x,y \in \bar{D}, \ x\neq y,
\end{equation}
 \begin{equation}\label{G_alpha_2}
\begin{aligned}
|G_{\alpha,v}(x,y,E)|\leq c_1(|x-y|^{2-d}),\ \  x,y\in\bar{D},\ {\rm for}\ \ d\ge 3,\\
|G_{\alpha,v}(x,y,E)|\leq c_1(|\ln\,|x-y||),\ \ x,y\in\bar{D},\ {\rm for}\ \ d=2,
\end{aligned}
\end{equation}
where $c_1=c_1(D,E,v,\alpha)>0.$
 
 Actually, properties \eqref{G_alpha_1}, \eqref{G_alpha_2} are well-known 
  for $\sin \alpha = 0$ (the case of the Direchlet boundary condition)
 and for $\cos \alpha = 0$ (the case of the Neumann boundary condition).
 Properties \eqref{G_alpha_1}, \eqref{G_alpha_2} with $d\geq 3$, $\frac{\sin \alpha}{\cos \alpha}<0$, $v\equiv 0$ and $E=0$ were proven in \cite{Lanzani2004}.
 For $d=2$ see also \cite{Begehr2010}.
 In Section 6 we give proofs of \eqref{G_alpha_1}, \eqref{G_alpha_2} for the case of general $\alpha$, $v$ and $E$.
 
 In addition, under assumptions of Theorem \ref{main}, the following identity holds 
 (see formula (3.9) of \cite{IN2012}):
\begin{equation}\label{Ident}
			\int\limits_{D} (v-v^0)\psi \psi^0 dx = \int\limits_{\partial D} [\psi]_\alpha  
		\left(\hat{M}_{\alpha,v} - \hat{M}_{\alpha,v^0}\right) [\psi^0]_\alpha  dx
\end{equation}
for all sufficiently regular solutions $\psi$, $\psi^0$ of equation \eqref{eq} in $\bar{D}$ for potentials $v$, $v^0$,
respectively, where $[\psi]_\alpha$, $[\psi^0]_\alpha $ are defined according to \eqref{[psi]_alpha}.

	Identity (\ref{Ident}) for  $\sin\alpha = 0$  is reduced to the Alessandrini identity (Lemma 1 of \cite{Alessandrini1988}).

We will use also that:
\begin{subequations}\label{R(k)}
\begin{equation}
	\begin{aligned}
	\|\hat{R}(k)u\|_{C^{1+\delta}(\Omega)} \leq c_2(D,\Omega,v,k,\delta)\|u\|_{\mathbb{L}^\infty(D)},\\
	\hat{R}(k)u(x) = \int\limits_D R(x,y,k) u(y) dy, \ \ \ x\in\Omega,\\ k \in \mathbb{C}^d\setminus(\mathbb{R}^d \cup \mathcal{E}),
	\end{aligned}
\end{equation}

\begin{equation}
	\begin{aligned}
	\|\hat{R}_\gamma(k)u\|_{C^{1+\delta}(\Omega)} \leq c_3(D,\Omega,v,k,\gamma,\delta)\|u\|_{\mathbb{L}^\infty(D)},\\
	\hat{R}_\gamma(k)u(x) = \int\limits_D R_\gamma(x,y,k) u(y) dy, \ \ \ x\in\Omega,\\ 
	\gamma\in \mathbb{S}^{d-1}, \ k \in \mathbb{R}^d\setminus(\{0\} \cup \mathcal{E}_\gamma),
	\end{aligned}
\end{equation}
\end{subequations}
for $u\in \mathbb{L}^\infty(D)$, $\delta\in [0,1)$, where  $\Omega$ is such an open bounded domain in $\mathbb{R}^d$ 
that $\bar{D}\in \Omega$ and
$C^{1+\delta}$ denotes $C^1$ with the first derivatives belonging to the H\"older space $C^{\delta}$.

	We will use also  the Green formula:
\begin{equation}\label{G_F}
		\int\limits_{\partial D} \left( \phi_1 \frac{\partial\phi_2}{\partial \nu} 
			- \phi_2 \frac{\partial\phi_1}{\partial \nu} \right)  dx = 
			\int\limits_{D} \left(\phi_1 \Delta \phi_2 - \phi_2 \Delta \phi_1\right) dx,
\end{equation}
	where $\phi_1$ and $\phi_2$ are arbitrary sufficiently regular functions in $\bar{D}$.


\subsection{Proof of Theorem \ref{main} and Proposition 3.1} 

For the case when $\sin \alpha = 0$, Theorem \ref{main} and Proposition 3.1 were proved in 
\cite{Novikov2005}. In this subsection we generalize the
proof of \cite{Novikov2005} to the case $\sin \alpha \neq 0$. 
We proceed from the following  formulas and equations (being valid under assumption \eqref{2.2} on $v^0$ and $v$):
\begin{equation}\label{delta_h}
	\begin{aligned}
		h(k,l) - 	h^0(k,l)=  \left(\frac{1}{2\pi}\right)^{d} 
		\int\limits_{\mathbb{R}^d} \psi^{0}(x,-l) (v(x) - v^0(x))\psi(x,k)dx,\\
  		k,l \in \mathbb{C}^d \setminus(\mathcal{E}^0\cup\mathcal{E}),\ k^2=l^2,\  |\mbox{Im}\,k| = |\mbox{Im} \,l| \neq 0,
 \end{aligned}
\end{equation}
\begin{equation}\label{psi_R}
	\begin{aligned}
			\psi(x,k) = \psi^0(x,k) + \int\limits_{\mathbb{R}^d} R^0(x,y,k) (v(y)-v^0(y))\psi(y,k)dy,\\
			x\in \mathbb{R}^d, \ k\in \mathbb{C}^d\setminus (\mathbb{R}^d \cup \mathcal{E}^0),
	\end{aligned}
\end{equation}
where \eqref{psi_R} at fixed $k$ is considered as an equation for $\psi = e^{ikx}\mu(x,k)$ with 
$\mu \in \mathbb{L}^{\infty}(\mathbb{R}^d)$;

\begin{equation}\label{delta_h^+}
	\begin{aligned}
		f(k,l) - 	f^0(k,l)=  \left(\frac{1}{2\pi}\right)^{d} 
		\int\limits_{\mathbb{R}^d} \psi^{+,0}(x,-l) (v(x) - v^0(x))\psi^+(x,k)dx,\\
  		k,l \in \mathbb{R}^d \setminus(\{0\}\cup \mathcal{E}^{+,0}\cup\mathcal{E}^+),\ k^2=l^2,
 \end{aligned}
\end{equation}
\begin{equation}\label{psi_R^+}
	\begin{aligned}
			\psi^+(x,k) = \psi^{+,0}(x,k) + \int\limits_{\mathbb{R}^d} R^{+,0}(x,y,k) (v(y)-v^0(y))\psi^+(y,k)dy,\\
			x\in \mathbb{R}^d, \ k\in \mathbb{R}^d\setminus (\{0\} \cup \mathcal{E}^{+,0}),
	\end{aligned}
\end{equation}
where \eqref{psi_R^+} at fixed  $k$	 is an equation for $\psi^+ \in \mathbb{L}^\infty(\mathbb{R}^d)$;

\begin{equation}\label{delta_h_gamma}
	\begin{aligned}
		h_\gamma(k,l)- h^0_\gamma(k,l) =  \left(\frac{1}{2\pi}\right)^{d} 
		\int\limits_{\mathbb{R}^d} \psi_{-\gamma}^{0}(x,-k,-l) (v(x) - v^0(x))\psi_\gamma(x,k)dx,\\
  		\gamma \in \mathbb{S}^{d-1},
  		k \in \mathbb{R}^d \setminus(\mathcal{E}_\gamma^0\cup\mathcal{E}_\gamma),\ l\in \mathbb{R}^d,\  k^2=l^2,  
 \end{aligned}
\end{equation}
\begin{equation}\label{psi_R_gamma}
	\begin{aligned}
			\psi_\gamma(x,k) = \psi_\gamma^0(x,k) + \int\limits_{\mathbb{R}^d} R^0_\gamma(x,y,k) (v(y)-v^0(y))\psi_\gamma(y,k)dy,\\
			x\in \mathbb{R}^d,\ \gamma \in \mathbb{S}^{d-1}, \ k\in \mathbb{R}^d\setminus (\{0\} \cup \mathcal{E}_\gamma^0),
	\end{aligned}
\end{equation}
where \eqref{psi_R_gamma} at fixed $\gamma$ and $k$ is considered as an equation for 
$\psi_\gamma \in \mathbb{L}^{\infty}(\mathbb{R}^d)$.

We recall that $\psi^+$, $f$, $\psi$, $h$, $\psi_\gamma$, $h_\gamma$ were defined in Sections 2, 3 by means of 
\eqref{2.3} - \eqref{2.9}, \eqref{psi(x,k,l)}.
Equation \eqref{psi_R^+} is well-known in the classical scattering theory for the Schr\"odinger equation \eqref{2.1}.
Formula \eqref{delta_h^+} was given, in particular, in \cite{Stefanov1990}. To our knowledge formula and equations
\eqref{delta_h},  \eqref{psi_R}, \eqref{psi_R_gamma} were given for the first time in \cite{Novikov1996}, 
whereas formula \eqref{delta_h_gamma}  was given for the first time in \cite{Novikov2005}.

In addition, under assumption \eqref{2.2} on $v^0$ and $v$:
\begin{subequations}\label{4.11}
	\begin{equation}
		\begin{aligned}
			\text{equation \eqref{psi_R} at fixed $k \in \mathbb{C}^d \setminus (\mathbb{R}^d\cup \mathcal{E}^0)$ 
			is uniquely solvable}\\
			\text{
			for $\psi = e^{ikx}\mu(x,k)$ with 
$\mu \in \mathbb{L}^{\infty}(\mathbb{R}^d)$ if and only if $k\notin \mathcal{E}$;}
		\end{aligned}
	\end{equation}
	\begin{equation}
		\begin{aligned}
			\text{equation \eqref{psi_R^+} at fixed $k \in \mathbb{R}^d \setminus (\{0\}\cup \mathcal{E}^{+,0})$ 
			is uniquely}\\
			\text{ solvable
			for $\psi^+ \in \mathbb{L}^{\infty}(\mathbb{R}^d)$ if and only if $k\notin \mathcal{E}^+$;}
		\end{aligned}
	\end{equation}
	\begin{equation}
		\begin{aligned}
			\text{equation \eqref{psi_R_gamma} at fixed $\gamma\in \mathbb{S}^{d-1}$ and 
			$k \in \mathbb{R}^d \setminus (\{0\}\cup \mathcal{E}^{+}_\gamma)$}\\
			\text{is uniquely solvable
			for $\psi_\gamma \in \mathbb{L}^{\infty}(\mathbb{R}^d)$ if and only if $k\notin \mathcal{E}_\gamma$.}
		\end{aligned}
	\end{equation}
\end{subequations}

Let us prove Theorem \ref{main} for the case of the Faddeev functions $\psi$, $h$. The proof of 
 Theorem \ref{main} for the case 
of $\psi^+$, $f$ and the proof of Proposition 3.1 
are similar.

Note that formula \eqref{eq_h} follows directly from \eqref{Ident} and \eqref{delta_h}.

Using \eqref{2.17} and applying \eqref{Ident} for equation \eqref{psi_R}, we get that
\begin{equation}\label{4.12}
	\begin{aligned}
		\psi(x,k) - \psi^0(x,k) = 
		\int\limits_{\partial D} \int\limits_{\partial D}
		[R^0(x,\xi,k)]_{\xi, \alpha}
		\left(M_{\alpha,v} - M_{\alpha,v^0}\right)(\xi,y,E) [\psi(y,k)]_\alpha d\xi dy, \\  x\in \mathbb{R}^d\setminus \bar{D},
	\end{aligned}
\end{equation}
where
\begin{equation}
	[R^0(x,\xi,k)]_{\xi, \alpha} = \left(\cos \alpha - \sin \alpha \frac{\partial}{\partial \nu_\xi}\right)R^0(x,\xi,k).
\end{equation}
	Equation (\ref{eq_psi_alpha}) follows from formula (\ref{4.12}), definition \eqref{[psi]_alpha} and the property that
	\begin{equation}
		\lim\limits_{\varepsilon \rightarrow +0} 
		\left(
		\cos \alpha - \sin \alpha \frac{\partial}{\partial \nu_x}
		\right)
		u(x+\varepsilon\nu_x) = [u(x)]_\alpha, \ \ \ x\in \partial D,
	\end{equation}
for $u(x) = \psi(x,k)-\psi^0(x,k)$.
\subsection{Proofs of Propositions 3.2 and 3.4}

In this subsection we prove Propositions 3.2, 3.4 for the case of equation \eqref{eq_psi_alpha} for
$[\psi]_\alpha$. The proofs of Propositions 3.2 and 3.4 for
the cases of $\psi^+$ and $\psi_\gamma$ are absolutely similar.
\vspace{3mm}

\begin{Proof} {\it Proposition 3.2.} 
 The proof of Proposition 3.2 for the case of $\sin \alpha =0$
was given in \cite{Novikov2005}. Let us assume that $\sin \alpha \neq 0$.

Using \eqref{4.2}, we find that
\begin{equation}\label{4.14}
	\left(M_{\alpha,v} - M_{\alpha,v^0}\right)(\xi,y,E) = 
	\frac{1}{\sin^2\alpha}\left(G_{\alpha,v} - G_{\alpha,v^0}\right)(\xi,y,E), \ \ \ \xi,y\in \partial D.
\end{equation}
Using \eqref{2.17}, \eqref{4.1}, \eqref{G_F} and the impedance boundary condition ({\ref{bound_cond}}) for $G_{\alpha,v}$, $G_{\alpha,v^0}$, 
we get that
\begin{equation}\label{4.15}
	\begin{aligned}
			\int\limits_{\partial D} [R^0(x,\xi,k)]_{\alpha,\xi} (G_{\alpha,v} - &G_{\alpha,v^0})(\xi,y,E) d\xi =\\
					= 
							\int\limits_{\partial D} \Bigg( [R^0(x,\xi,k)&]_{\alpha,\xi}\left(G_{\alpha,v} - G_{\alpha,v^0}\right)(\xi,y,E) d\xi - \\ 
							&-
			 	R^0(x,\xi,k) [\left(G_{\alpha,v} - G_{\alpha,v^0}\right)(\xi,y,E)]_{\alpha,\xi} \Bigg) d\xi =\\
			 	=
			 	\sin \alpha \int\limits_{D} \Bigg( R^0(x,&\xi,k) \Delta_\xi \left(G_{\alpha,v} - G_{\alpha,v^0}\right)(\xi,y,E) d\xi -\\
			 	&-
			 	\left(G_{\alpha,v} - G_{\alpha,v^0}\right)(\xi,y,E) \Delta_\xi R^0(x,\xi,k)  \Bigg) d\xi = \\
			 	= 	\sin \alpha \int\limits_{D} R^0(x,\xi,k) &\left(v(\xi) - v^0(\xi)\right) G_{\alpha,v} (\xi,y,E) d\xi, \ \ \
			 	x\in \mathbb{R}^d \setminus \bar{D}, \ y \in \partial D. 
	\end{aligned}
\end{equation}
Combining \eqref{4.12}, \eqref{4.14} and \eqref{4.15}, we obtain that
\begin{equation}\label{4.16}
		A_\alpha(x,y,k) = \lim_{\varepsilon \rightarrow +0}
		\left(\cos \alpha - \sin \alpha \frac{\partial}{\partial \nu_x}\right)
		B_\alpha(x+\varepsilon\nu_x,y,k), \ \ \ x,y\in \partial D,		
\end{equation}
where
\begin{equation}\label{4.17}
	\begin{aligned}
	B_\alpha(x,y,k) = \int\limits_{\partial D} [R^0(x,\xi,k)]_{\alpha,\xi} (M_{\alpha,v} - M_{\alpha,v^0})(\xi,y,E) d\xi=\\
	=\frac{1}{\sin\alpha} \int\limits_{D} R^0(x,\xi,k)\left(v(\xi)-v^0(\xi)\right) G_{\alpha,v}(\xi,y,E)  d\xi, \\
	x\in \mathbb{R}^d\setminus \bar{D},\ y\in \partial D.
	\end{aligned}
\end{equation}
Thus, we have that the limit in \eqref{4.16} (and, hence, in (\ref{A_alpha})) is well defined and 
\begin{equation}\label{4.18}
		A_\alpha(x,y,k) = \frac{1}{\sin\alpha} \int\limits_{D} [R^0(x,\xi,k)]_{x,\alpha}
		\left(v(\xi)-v^0(\xi)\right) G_{\alpha,v}(\xi,y,E)  d\xi, \ \
	x,y\in \partial D.
\end{equation}
Let $\hat{A}_\alpha(k)$ denote the linear integral operator on $\partial D$ with the Schwartz kernel
$A_{\alpha}(x,y,k)$ of \eqref{A_alpha}, \eqref{4.18}. Using \eqref{G_alpha_2}, \eqref{R(k)}, \eqref{4.18},
we obtain that
\begin{equation}\label{A(k)}
	\begin{aligned}
	\hat{A}_\alpha(k):\ \mathbb{L}^\infty(\partial D) \rightarrow C^\delta(\partial D) \\
		\text{is a bounded linear operator.}
	\end{aligned}
\end{equation}
As a corollary of \eqref{A(k)}, $\hat{A}_\alpha(k)$ is a compact operator in $\mathbb{L}^\infty(D)$.
\end{Proof}

\begin{Proof}{ \it Proposition 3.4.}
Using \eqref{2.17}, \eqref{eq_phi} and \eqref{G_F}, we get that
\begin{equation}\label{4.19}
\begin{aligned}
	\int\limits_{\partial D}  
	\left( 
		\phi_{\alpha}(\xi,y) \frac{\partial}{\partial \nu_\xi}R^0(x,\xi,k)   -
		 R^0(x,\xi,k) \frac{\partial}{\partial \nu_\xi}  \phi_{\alpha}(\xi,y) 
	\right)d\xi =\\
=
\int\limits_{D}  
	\left( 
		\phi_{\alpha}(\xi,y)  \Delta_\xi R^0(x,\xi,k)   -
		 R^0(x,\xi,k) \Delta_\xi  \phi_{\alpha}(\xi,y) 
	\right)d\xi =\\=
\int\limits_{D}  R^0(x,\xi,k) (v^0(\xi)-E+\lambda) \phi_{\alpha}(\xi,y) d\xi, \\  x\in \mathbb{R}^d \setminus \bar{D}, \ 
y\in \partial D.
	\end{aligned}
	\end{equation}
 Combining \eqref{eq_phi}, \eqref{4.17}
 and \eqref{4.19}, we find that 
 \begin{equation}\label{4.20}
 	\begin{aligned}
		B_\alpha(x,y,k) = 	\int\limits_{\partial D} [R^0(x,\xi,k)]_{\xi,\alpha} \phi_{\alpha}(\xi,y) d\xi &=\\ = 
		\int\limits_{\partial D} R^0(x,\xi,k) [\phi_{\alpha}(\xi,y)]_{\xi,\alpha} d\xi 
		- \sin \alpha \int\limits_{D}  &R^0(x,\xi,k) (v^0(\xi)-E+\lambda) \phi_{\alpha}(\xi,y) d\xi, \\
		&\ \ \ \ \text{  } \ \ \ \text{  } \ \ \  x\in \mathbb{R}^d\setminus \bar{D},\ y\in \partial D.
	\end{aligned}
 \end{equation}
 Combining  \eqref{4.16} and \eqref{4.20}, we obtain \eqref{3.8}.
 
Formula \eqref{3.10} follows from \eqref{eq_phi} and the definition of $\hat{\Phi}$.

\end{Proof}
\section{Proof of Proposition 3.3}
For the case when $\sin \alpha = 0$, Proposition 3.3 was proved in 
\cite{Novikov2005}. In this section we prove Proposition 3.3
for $\sin \alpha \neq 0$. 
We will prove  Proposition 3.3  
for the case of equation \eqref{eq_psi_alpha} for $[\psi]_\alpha$. The proofs  for
the cases of $\psi^+$ and $\psi_\gamma$ are similar.

According to \eqref{4.11}, to prove Proposition 3.3 (for the case of $\psi$) it is sufficient to show
that 
 equation  (\ref{eq_psi_alpha}) (at fixed $k\in \mathbb{C}^d\setminus (\mathbb{R}^d\cup \mathcal{E}^0)$)  is  uniquely
solvable in the space of bounded functions on  $\partial D$ if and only if equation \eqref{psi_R}
	is uniquely solvable for $\psi = e^{ikx}\mu(x,k)$ with  $\mu\in \mathbb{L}^\infty(\mathbb{R}^d)$.
	
Let equation (\ref{psi_R}) have several solutions. Then, repeating the proof of Theorem \ref{main}  separately 
for each solution, we find that $[\psi]_\alpha$ on $\partial D$ for each of these solutions
satisfies equation (\ref{eq_psi_alpha}). Thus, using also \eqref{correct_M} we obtain that 
equation (\ref{eq_psi_alpha}) has at least as many solutions as equation
\eqref{psi_R}.

To prove the converse (and thereby to prove Proposition 3.3) it remains to show that any solution $[\psi]_\alpha$
of (\ref{eq_psi_alpha}) can be continued to a continuos solution of \eqref{psi_R}.

Let $\psi$ be the solution of \eqref{eq} with the impedance boundary data $[\psi]_\alpha$, satisfying
\eqref{eq_psi_alpha}. Let
\begin{equation}\label{5.1}
  \psi_1(x) 
  =    \psi^0(x,k) + \int\limits_{D} R^0(x,y,k) (v(y)- v^0(y)) \psi(y) dy, \ \ \ x\in \mathbb{R}^d.
\end{equation}
Using \eqref{R(k)}, 
we obtain that 
\begin{equation}\label{5.2+}
\text{$\psi_1$ defined by \eqref{5.1} belongs to $C^{1+\delta}(\mathbb{R}^d)$, $\delta \in [0,1)$.} 
\end{equation}
We have that
\begin{equation}\label{5.2}
	(-\Delta + v^0(x) - E)\psi(x) = (v^0(x)- v(x)) \psi(x), \ \  \ x\in D,
\end{equation}
\begin{equation}\label{5.3}
	\begin{aligned}
	(-\Delta + v^0(x) - E)\psi_1(x) 
		= \int\limits_{D} -\delta(x-y)  
		(v(y)- v^0(y)) \psi(y) dy =\\= (v^0(x)- v(x)) \psi(x), \ \ \ x\in D.
	\end{aligned}
\end{equation}
Combining \eqref{Ident} and \eqref{4.17}, we get that
\begin{equation}\label{5.10}
	\int\limits_{D} R^0(x,y,k) (v(y)- v^0(y)) \psi(y) dy 
	 = \int\limits_{\partial D} B_\alpha(x,y,k) [\psi(y)]_\alpha dy, \ \ \ x\in \mathbb{R}^d\setminus \bar{D}.
\end{equation}
Using \eqref{eq_psi_alpha}, \eqref{4.16}, \eqref{5.2+}, \eqref{5.10}, we find that
\begin{equation}\label{psi1psi}
	[\psi_1(x)]_\alpha = [\psi^0(x,k)]_\alpha + \int\limits_{\partial D} A_\alpha(x,y,k) [\psi(y)]_\alpha dy =[\psi(x)]_\alpha, 
	\ \ \ x\in \partial D.
\end{equation}
Using \eqref{5.2}, \eqref{5.3} and \eqref{psi1psi}, we obtain that
\begin{equation}
	\begin{aligned}
	(-\Delta + v^0(x) - E) (\psi_1(x)-\psi(x)) = 0, \ \ \ x\in D,\\
		[\psi_1(x)-\psi(x)]_\alpha = 0, \ \ \ x\in\partial D.
	\end{aligned}
\end{equation}
Since $v^0$ satisfies (\ref{correct_M}), we get that 
\begin{equation}\label{psi=psi1}
\psi_1(x) = \psi(x), \ \ \ x\in \bar{D}.
\end{equation} 
Combining \eqref{5.1}, \eqref{5.2+} and \eqref{psi=psi1}, we find that  $\psi_1$ is a continuos solution of \eqref{psi_R}.

\section{Proofs of properties (\ref{G_alpha_1}), (\ref{G_alpha_2})}

As it was mentioned in Subsection 4.1,  properties (\ref{G_alpha_1}), (\ref{G_alpha_2}) are well-known for $\cos \alpha =0$
(the case of the Neumann boundary condition). To extend these properties to the case
of general $\alpha$, $v$, $E$, we use the following schema:
\begin{enumerate}
\item $G_{\alpha_1,v} \rightarrow G_{\alpha_2,v}$ by means of Lemma 6.1 given bellow
(with $\sin\alpha_1\neq 0$ and $\sin\alpha_2\neq 0$),
\item $G_{\alpha,v_1} \rightarrow G_{\alpha,v_2}$ by means of Lemma 6.2 given bellow.
\end{enumerate}
The proofs of steps 1, 2 are based on the theory of Fredholm linear integral equations of the second kind. 

Starting from (\ref{G_alpha_1}), (\ref{G_alpha_2}) for $\cos \alpha =0$ and combining steps 1, 2 and
the property
\begin{equation}
	G_{\alpha,v}(\cdot,\cdot,E) = G_{\alpha,v-E}(\cdot,\cdot,0),
\end{equation}  
we obtain  these properties for the case when $\sin \alpha \neq 0$.

As it was already mentioned in Section 4, properties (\ref{G_alpha_1}), (\ref{G_alpha_2}) are well-known for $\sin \alpha = 0$ 
(the case of the Dirichlet boundary condition).

\begin{Lemma}\label{Lemma6.1}
	Let $D$ satisfy \eqref{eq_c} and potential $v$ satisfy (\ref{eq_c2}), (\ref{correct_M}) 
	for some fixed $E$  and  for $\alpha = \alpha_1$, $\alpha = \alpha_2$ simultaneously, where $\sin \alpha_1 \neq 0$
	and $\sin \alpha_2 \neq 0$. 
	Let $G_j$ denote the Green function $G_{\alpha_j,v}$, $j=1,2$. Let $G_1$ satisfy:
	\begin{equation}\label{G_1c}
			G_1(x,y,E)\ \ {\rm is\ continuous \ in }\ \ x,y \in \bar{D},\  x\neq y,
	\end{equation}
 	\begin{equation}\label{G_1e}
			\begin{aligned}
				|G_1(x,y,E)| \leq a_1|x-y|^{2-d}\ \ \ {\rm for}\ \ d\ge 3,\\
				|G_1(x,y,E)|\leq a_1|\ln\,|x-y||\ \ \ {\rm for}\ \ d=2,\\
				x,y\in \bar{D}.
			\end{aligned}
	\end{equation}
	Then:
		\begin{equation}\label{G_2c}
			G_2(x,y,E)\ \ {\rm is\ continuous \ in }\ \ x,y \in \bar{D},\  x\neq y,
	\end{equation}
 	\begin{equation}
			\begin{aligned}\label{G_2e}
				|G_2(x,y,E)| \leq a_2|x-y|^{2-d}\ \ \ {\rm for}\ \ d\ge 3,\\
				|G_2(x,y,E)|\leq a_2|\ln\,|x-y||\ \ \ {\rm for}\ \ d=2,\\
				x,y\in \bar{D},
			\end{aligned}
	\end{equation}
	where $a_2 =a_2(D,E,a_1,v, \alpha_1,\alpha_2)>0$. 
\end{Lemma}

\begin{Proof} {\it Lemma \ref{Lemma6.1}.}
First, we derive formally some formulas and equations relating the Green functions $G_1$ and $G_2$. Then, proceeding 
from these formulas and equations, we obtain, in particular, estimates \eqref{G_2c}, \eqref{G_2e}.

Consider $W = G_2 - G_1$. Using definitions of $G_1$, $G_2$ and formula \eqref{psiG+++}, we find that:
	\begin{equation}\label{eq_W}
		(-\Delta_x+v(x) - E) W(x,y) = 0, \ \ x,y\in D,
	\end{equation} 
	\begin{equation}\label{bc_W}
		\begin{aligned}
		\bigg(\cos \alpha_2 \,W(x,y) - \sin\alpha_2\, \frac{\partial W}{\partial \nu_x}(x,y)\bigg) \Big|_{x\in \partial D} &= \\
		= -\bigg(\cos \alpha_2 \,G_1(x,y,E) - \sin\alpha_2\, &\frac{\partial G_1}{\partial \nu_x}(x,y,E)\bigg) 
		\Big|_{x\in \partial D} =\\
		= -\bigg(\cos \alpha_2 \,G_1(x,y,E) - \sin\alpha_2\, &\frac{\cos\alpha_1}{\sin\alpha_1}G_1(x,y,E)\bigg)\Big|_{x\in \partial D}=\\
		=&\frac{\sin(\alpha_2-\alpha_1)}{\sin\alpha_1}G_1(x,y,E)\Big|_{x\in \partial D}, \ \ \ y\in D,
		\end{aligned}
	\end{equation}
	\begin{equation}\label{6.9}
		\begin{aligned}
		W(x,y) =  \frac{1}{\sin \alpha_1}\int\limits_{\partial D} 
	 \Big( \cos \alpha_1 W(\xi,y) - \sin \alpha_1 \frac{\partial W}{\partial \nu_\xi} (\xi,y)
	 \Big)G_1(\xi,x,E)  d\xi, \ \ \ &x,y\in D.
		\end{aligned}
	\end{equation}
 Using	\eqref{bc_W} and \eqref{6.9}, we find the following linear integral equation
 for $W(\cdot, y)$ on $\partial D$:
	\begin{equation}\label{eq_K1}
		W(\cdot,y) = W_0(\cdot,y) +   \hat{K}_1 W(\cdot, y) , \ \ \  y\in D,
	\end{equation}
	where
	\begin{equation}\label{def_W0}
		W_0(x,y) = \frac{\sin(\alpha_2-\alpha_1)}{\sin{\alpha_2}} \int\limits_{\partial D} 
		G_1(\xi,x,E) G_1(\xi,y,E) d \xi,
	\end{equation}
	\begin{equation}
		\begin{aligned}
		\hat{K}_1 u (x) = \frac{\sin(\alpha_2-\alpha_1)}{\sin{\alpha_2}\sin{\alpha_1}} \int\limits_{\partial D} 
		 G_1(\xi,x,E) u(\xi) d\xi, \\  x\in \partial D, \ y\in D, \text{ $u$ is a test function.}
		 \end{aligned}
	\end{equation}
In addition, for 
\begin{equation}\label{deltW}
	\delta_n W = W - \sum\limits_{j=1}\limits^{n} (\hat{K}_1)^{j-1} W_0
\end{equation}
equation \eqref{eq_K1} takes the form
	\begin{equation}\label{eq_deltaW}
		\delta_n W = (\hat{K}_1)^n W_0 + \hat{K}_1 \delta_n W.
	\end{equation}

 Our analysis based on \eqref{eq_W}-\eqref{eq_deltaW}
is given bellow.

		Using \eqref{G_1c}, \eqref{G_1e}, we obtain that
	\begin{equation}\label{suff_n}
		(\hat{K}_1)^n W_0 \in C(\partial D \times \bar{D}) \text{ for sufficiently great $n$ with respect to $d$},
	\end{equation}
	\begin{equation}
		\text{ $\hat{K}_1$ is a compact operator in $C(\partial D)$.}
	\end{equation}

Let us show that the homogeneous equation
\begin{equation}\label{6.13}
	u = \hat{K}_1 u, \ \ u\in C(\partial D),
\end{equation}
has only trivial solution $u\equiv 0$. 

Using the fact that the potential $v$ satisfy  (\ref{correct_M}) 
  for $\alpha = \alpha_1$,
we define $\psi$   by
\begin{equation}
	\begin{aligned}
		(-\Delta   + v(x) - E)\psi(x) = 0, \ \ \  x \in  D,\\
		\cos\alpha_1 \psi|_{\partial D}  - \sin \alpha_1 \frac{\partial\psi}{\partial \nu}|_{\partial D} = u.
	\end{aligned}
\end{equation}
Due to \eqref{psiG+++}, we have that
\begin{equation}\label{6.15}
	\psi(x) = \frac{1}{\sin\alpha_1} \int\limits_{\partial D} 
	 (\cos\alpha_1 \psi(\xi) - \sin\alpha_1 \frac{\partial \psi}{\partial \nu}(\xi)) G_1(\xi,x,E)d\xi, \ \ \ x\in D.
\end{equation}
Using \eqref{6.13}, \eqref{6.15}, we find that
\begin{equation}
	\frac{\sin(\alpha_2-\alpha_1)}{\sin\alpha_2} \psi(x) = \hat{K}_1 u (x) = u(x), \ \ \ x\in \partial D.
\end{equation}
Therefore, we have that
\begin{equation}\label{6.17}
	\cos\alpha_1 \psi(x)  - \sin \alpha_1 \frac{\partial\psi}{\partial \nu}(x) = 
	\frac{\sin(\alpha_2-\alpha_1)}{\sin\alpha_2} \psi(x), \ \ \ x\in \partial D.
\end{equation}
Since $\sin\alpha_1 \neq 0$ and $\sin\alpha_2 \neq 0$, using \eqref{6.17}, we obtain that
\begin{equation}
	\cos\alpha_2 \psi(x)  - \sin \alpha_2 \frac{\partial\psi}{\partial \nu}(x) = 0
\end{equation}
 Taking into account the fact that the potential $v$ satisfy  (\ref{correct_M}) 
  for $\alpha = \alpha_2$, we get that $\psi\equiv 0$ and $u \equiv 0$.   
  
  Proceeding from 
  \begin{equation}
		\begin{aligned}  
    F = W(x,y) \ \text{ and }\ 
  	F' = \frac{\cos\alpha_2}{\sin\alpha_2} W(x,y) -\frac{\sin(\alpha_2 -\alpha_1)}{\sin\alpha_1 \sin\alpha_2} G_1(x,y,E), 
  	\\ x \in \partial D,\ y\in \bar{D},
  	\end{aligned}
  \end{equation}
  found from \eqref{eq_K1}, \eqref{eq_deltaW} and \eqref{bc_W} (with $F'$ substituted in place of $\partial W/\partial \nu_x$), we consider 
  \begin{equation}
  	W(x,y) = \frac{1}{\sin \alpha_1}\int\limits_{\partial D} 
	 \Big( \cos \alpha_1 F(\xi,y) - \sin \alpha_1 F'(\xi,y)
	 \Big)G_1(\xi,x,E)  d\xi, \ \ \ x,y\in \bar{D}.
  \end{equation}
  Using \eqref{eq_K1} and properties of $G_1$ (including formula \eqref{psiG+++}), we subsequently obtain that
    \begin{equation}\label{6.24}
  	\lim\limits_{\varepsilon \rightarrow +0}
  	W(x - \varepsilon\nu_x,y) = F(x,y), \ \ \ x\in \partial D,\ y\in \bar{D}, 
  \end{equation}
  \begin{equation}
  	 \text{ $W$ satisfies \eqref{eq_W}, }
  \end{equation}
  \begin{equation}\label{6.26}
  	\lim\limits_{\varepsilon \rightarrow +0}
  	\frac{\partial}{\partial \nu_x}W(x - \varepsilon\nu_x,y) = F'(x,y),
  	\ \ \ x\in \partial D,\ y\in \bar{D}. 
  \end{equation}

  From \eqref{G_1c}, \eqref{G_1e}, \eqref{def_W0}-\eqref{6.13}, \eqref{6.24}-\eqref{6.26} it follows that $G_2$ defined as $G_2 = G_1+W$
  is the Green function for the operator $\Delta-v+E$ in $D$ with the impedance boundary condition \eqref{bound_cond} for $\alpha= \alpha_2$
  and that $G_2$ satisfies \eqref{G_2c}, \eqref{G_2e}.
  \end{Proof}
\begin{Lemma}\label{Lemma6.2}
	Let $D$ satisfy \eqref{eq_c} and potentials $v_1$, $v_2$ satisfy (\ref{eq_c2}), (\ref{correct_M}) 
	for some fixed $E$  and  $\alpha$. 
	Let $G_j$ denote the Green function $G_{\alpha,v_j}$, $j=1,2$. Let $G_1$ satisfy:
	\begin{equation}\label{G_1c2}
			G_1(x,y,E)\ \ {\rm is\ continuous \ in }\ \ x,y \in \bar{D},\  x\neq y,
	\end{equation}
 	\begin{equation}\label{G_1e2}
			\begin{aligned}
				|G_1(x,y,E)| \leq a_3|x-y|^{2-d}\ \ \ {\rm for}\ \ d\ge 3,\\
				|G_1(x,y,E)|\leq a_3|\ln\,|x-y||\ \ \ {\rm for}\ \ d=2,\\
				x,y\in \bar{D}.
			\end{aligned}
	\end{equation}
	Then:
		\begin{equation}\label{G_2c2}
			G_2(x,y,E)\ \ {\rm is\ continuous \ in }\ \ x,y \in \bar{D},\  x\neq y,
	\end{equation}
 	\begin{equation}\label{G_2e2}
			\begin{aligned}
				|G_2(x,y,E)| \leq a_4|x-y|^{2-d}\ \ \ {\rm for}\ \ d\ge 3,\\
				|G_2(x,y,E)|\leq a_4|\ln\,|x-y||\ \ \ {\rm for}\ \ d=2,\\
				x,y\in \bar{D},
			\end{aligned}
	\end{equation}
	where $a_4 =a_4(D,E,a_3,v_1, v_2, \alpha)>0$. 
\end{Lemma}

\begin{Proof} {\it Lemma \ref{Lemma6.2}.}
First, we derive formally some formulas and equations relating the Green functions $G_1$ and $G_2$. Then, proceeding 
from these formulas and equations, we obtain, in particular, estimates \eqref{G_2c2}, \eqref{G_2e2}.

Using (\ref{4.1}), the impedance boundary condition for $G_1$, $G_2$, we find that
	\begin{equation}\label{6.2.1}
		\begin{aligned}
			G_1(x,y,E) = \int\limits_{D} 
				 G_1(x,\xi,E) \Big( \Delta_\xi - v_2(\xi) + E\Big) G_2(\xi,y,E)\, d\xi,	\\
			G_2(x,y,E) = \int\limits_{D} 
				G_2(\xi,y,E) \Big(\Delta_\xi - v_1(\xi) + E\Big) G_1(x,\xi,E) \, d\xi, \\
		\int\limits_{\partial D} \left( G_1(x,\xi,E)  \frac{\partial G_2}{\partial \nu_\xi} (\xi,y,E)
		- G_2(\xi,y,E) \frac{\partial G_1}{\partial \nu_\xi}(x,\xi,E)  \right)  d\xi = 0,\\
		x, y \in D. 
		\end{aligned}
	\end{equation}
	
	Combining (\ref{6.2.1}) with (\ref{G_F}), we get that
	\begin{equation}\label{6.2.2}
		\begin{aligned}
			G_2 (\cdot,y,E)  -  G_1(\cdot,y,E) =
											 \hat{K}_2 G_2(\cdot,y,E),\ \ \ y \in D, 
		\end{aligned}
	\end{equation}
	where
	\begin{equation}\label{6.33}
		\hat{K}_2 u(x) = \int\limits_{D} \left( v_2(\xi) -v_1(\xi) \right) G_1(x,\xi,E) u(\xi) d\xi.
	\end{equation}
	In addition, for 
\begin{equation}\label{deltG}
	\delta_n G = G_2 - \sum\limits_{j=1}\limits^{n} (\hat{K}_2)^{j-1} G_1
\end{equation}
equation \eqref{6.2.2} takes the form
	\begin{equation}\label{eq_deltG}
		\delta_n G = (\hat{K}_2)^n G_1 + \hat{K}_2 \delta_n G.
	\end{equation}
	
	 Our analysis based on \eqref{6.2.1}-\eqref{eq_deltG}
is given bellow.
	
	Using \eqref{G_1c2}, \eqref{G_1e2}, we find that
	\begin{equation}\label{6.37}
		(\hat{K}_2)^n G_1 \in C(\bar D \times \bar{D}) \text{ for sufficiently great $n$ with respect to $d$},
	\end{equation}
	\begin{equation}\label{6.38}
		\text{ $\hat{K}_2$ is a compact operator in $C(\bar{D})$.}
	\end{equation}
	
	Let us show that the homogeneous equation
\begin{equation}\label{6.38a}
	u = \hat{K}_2 u, \ \ u\in C(\bar{D}),
\end{equation}
has only trivial solution $u \equiv 0$.
Using \eqref{6.33}, \eqref{6.38a} and properties of the Green function $G_1$, we find that
\begin{equation}
	\begin{aligned}
	(-\Delta+v_1(x) - E) u(x) = \int\limits_D - \delta(x-\xi)\left( v_2(\xi) -v_1(\xi) \right)  u(\xi) d\xi =\\
		= (v_1 - v_2) u(x), \ \ \ x\in D,\\
		\cos\alpha\, u(x)  - \sin \alpha \frac{\partial u}{\partial \nu}(x) = 0, \ \ \ x\in \partial D.
	\end{aligned}
\end{equation}
	Using \eqref{G_1c2}, \eqref{G_1e2}, we find that $u \in C(\bar{D})$. Taking into account the fact that the potential $v_2$ satisfy  (\ref{correct_M}), we get that $u \equiv 0$. 
	
	Proceeding from \eqref{G_1c2}, \eqref{G_1e2}, \eqref{6.37}, \eqref{6.38} it follows that $G_2$ found from \eqref{6.2.2}, \eqref{eq_deltG}
	is the Green function for the operator $\Delta-v+E$ in $D$ with the impedance boundary condition \eqref{bound_cond} for $v = v_2$ and that
	$G_2$ satisfies \eqref{G_2c2}, \eqref{G_2e2}.
\end{Proof}

\section*{Acknowledgements}
This work was partially supported by TFP No 14.A18.21.0866 of Ministry of Education
and Science of Russian Federation. The second author was also partially supported by the Russian Federation Goverment
grant No. 2010-220-01-07

\noindent
{ {\bf M.I. Isaev}\\
Centre de Math\'ematiques Appliqu\'ees, Ecole Polytechnique,

91128 Palaiseau, France\\
Moscow Institute of Physics and Technology,

141700 Dolgoprudny, Russia\\
e-mail: \tt{isaev.m.i@gmail.com}}\\

\noindent
{ {\bf R.G. Novikov}\\
Centre de Math\'ematiques Appliqu\'ees, Ecole Polytechnique,

91128 Palaiseau, France\\
Institute of Earthquake Prediction Theory and Mathematical Geophysics RAS,

117997 Moscow, Russia\\ 
e-mail: \tt{novikov@cmap.polytechnique.fr}}

\end{document}